\documentclass[11pt]{article}
\usepackage{cite}
\usepackage{mathrsfs}
\usepackage{amsfonts}
\usepackage{amsmath}
\usepackage{amsfonts,amssymb,color}
\usepackage{dsfont}
\usepackage{curves}
\usepackage{mathrsfs}
\usepackage{pifont}
\usepackage{amssymb}
\usepackage{latexsym,amsmath,amssymb,amsfonts,epsfig,graphicx,cite,psfrag}
\usepackage{eepic,color,colordvi,amscd}

\newtheorem{theorem}{Theorem}[section]

\newtheorem{lemma}{Lemma}[section]

\newtheorem{conjecture}{Conjecture}[section]

\newcommand{\qed}{\hfill\rule{0.5em}{0.809em}}

\def\emptyset{\mbox{{\rm \O}}}
\def\bar{\overline}

\def\qed{\hfill \rule{4pt}{7pt}}

\def\pf{\noindent {\it Proof. }}

\begin{document}
\title{Optimal chromatic bound for ($P_3\cup P_2$, house)-free graphs}

\author{{Rui Li $^{a,*}$\footnote{Email address: lirui@hhu.edu.cn}\quad  Di Wu$^{b,}$\footnote{Email address: 1975335772@qq.com
}\quad Jinfeng Li $^{a,}$\footnote{Email address: 1345770246@qq.com}}\\
       {\small $^a$ School of Mathematics, Hohai University}\\
        {\small 8 West Focheng Road, Nanjing, 211100, China}\\
        {\small $^b$ School of Mathematical Science, Nanjing Normal University}\\
        {\small 1 Wenyuan Road, Nanjing, 210046, China}
        }
\date{}
\maketitle

\begin{abstract}
Let $G$ and $H$ be two vertex disjoint graphs. The {\em union} $G\cup H$ is the graph with $V(G\cup H)=V(G)\cup V(H)$ and $E(G\cup H)=E(G)\cup E(H)$. We use $P_k$ to denote a {\em path} on $k$ vertices, use {\em house} to denote the complement of $P_5$.  In this paper, we show that $\chi(G)\le2\omega(G)$ if $G$ is ($P_3\cup P_2$, house)-free. Moreover, this bound is optimal when $\omega(G)\ge2$.

\begin{flushleft} {\em Key words and
phrases:}  chromatic number; clique number; $\chi$-binding function; $(P_3\cup P_2)$-free graphs 

\end{flushleft}
\end{abstract}

\section{Introduction}

All graphs considered in this paper are finite and simple. We use $P_k$ and $C_k$ to denote a {\em path} and a {\em cycle} on $k$ vertices respectively, and follow \cite{BM1976} for undefined notations and terminology. Let $G$ be a graph, and $X$ be a subset of $V(G)$. We use $G[X]$ to denote the subgraph of $G$ induced by $X$, and call $X$ a {\em clique} ({\em independent set}) if $G[X]$ is a complete graph (has no edge). The {\em clique number} $\omega(G)$ of $G$ is the maximum size taken over all cliques of $G$.

For $v\in V(G)$, let $N_G(v)$ be the set of vertices adjacent to $v$, $d_G(v)=|N_G(v)|$, $N_G[v]=N_G(v)\cup \{v\}$, $M_G(v)=V(G)\setminus N_G[v]$. For $X\subseteq V(G)$, let $N_G(X)=\{u\in V(G)\setminus X\;|\; u$ has a neighbor in $X\}$ and $M_G(X)=V(G)\setminus (X\cup N_G(X))$. If it does not cause any confusion, we will omit the subscript $G$ and simply write $N(v), d(v), N[v], M(v), N(X)$ and $M(X)$.  Let $\Delta(G)$ ($\delta(G)$) denote the maximum (minimum) degree of $G$.

Let $G$ and $H$ be two vertex disjoint graphs. The {\em union} $G\cup H$ is the graph with $V(G\cup H)=V(G)\cup (H)$ and $E(G\cup H)=E(G)\cup E(H)$. The union of $k$ copies of the same graph $G$ will be denoted by $kG$. The {\em join} $G+H$ is the graph with $V(G+H)=V(G)+V(H)$ and $E(G+H)=E(G)\cup E(H)\cup\{xy\;|\; x\in V(G), y\in V(H)$$\}$. The complement of a graph $G$ will be denoted by $\bar G$.  We say that $G$ induces $H$ if $G$ has an induced subgraph isomorphic to $H$, and say that $G$ is $H$-free otherwise. Analogously, for a family $\cal H$ of graphs, we say that $G$ is ${\cal H}$-free if $G$ induces no member of ${\cal H}$.

Let $k$ be a positive integer, and let $[k]=\{1, 2, \ldots, k\}$. A $k$-{\em coloring} of $G$ is a mapping $c: V(G)\mapsto [k]$ such that $c(u)\neq c(v)$ whenever $u\sim v$ in $G$. The {\em chromatic number} $\chi(G)$ of $G$ is the minimum integer $k$ such that $G$ admits a $k$-coloring.  It is certain that $\chi(G)\ge \omega(G)$. A {\em perfect graph} is one such that  $\chi(H)=\omega(H)$ for all of its induced subgraphs $H$.  A family $\cal G$ of graphs is said to be $\chi$-{\em bounded} if there is a function $f$ such that $\chi(G)\le f(\omega(G))$ for every $G\in\cal G$, and if such a function does exist for $\cal G$, then $f$ is said to be a {\em binding function} of $\cal G$ \cite{Gy75}.

Erd\"{o}s \cite{ER59} proved that for any positive integers $k,l\ge3$, there exists a graph $G$ with $\chi(G)\ge k$ and no cycles of length less than $l$. This result motivates us to study the chromatic number of $F$-free graphs, where $F$ is a forest (a disjoint union of trees). Gy\'{a}rf\'{a}s \cite{Gy75} and Sumner \cite{Su81} independently, proposed the following famous conjecture.

\begin{conjecture}\label{tree}\cite{Gy75, Su81}
	Let $F$ be a forest. Then $F$-free graphs are $\chi$-bounded.
\end{conjecture}

Since $P_4$-free graphs are perfect, the class of $P_5$-free graphs has attracted a great deal of interest in recent years.  Up to now, the best known $\chi$-binding function for $P_5$-free graphs is $f(\omega)=\omega(G)^{\log_2\omega(G)}$ \cite{SS21}. We refer the interested readers to \cite{GM22} for results of $P_5$-free graphs, and to \cite{RS04,SR19,SS20} for more results and problems about the $\chi$-bounded problem. Notice that $P_5$ is a connected graph with five vertices. In Conjecture \ref{tree}, $F$ is a forest, what if we consider a disconnected graph with five vertices? Therefore, $(P_3\cup P_2)$-free graphs aroused our interest. The best known $\chi$-binding function for $(P_3\cup P_2)$-free graphs is $f(\omega)=\frac{1}{6}\omega(\omega+1)(\omega+2)$ \cite{BA18}.  Although the function is polynomial, we want to find a better $\chi$-binding function, preferably optimal, for some subclass of $(P_3\cup P_2)$-free graphs.

\begin{figure}[htbp]\label{fig-1}
	\begin{center}
		\includegraphics[width=15cm]{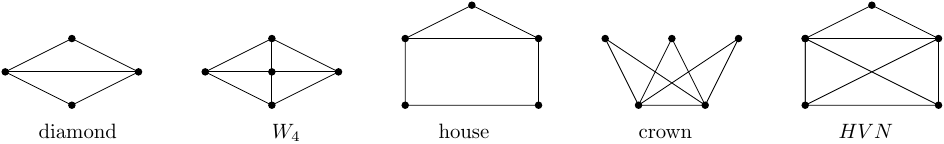}
	\end{center}
	\vskip -15pt
	\caption{Illustration of diamond, $W_4$, house, crown, and $HVN$.}
\end{figure}

A {\em diamond} is the graph $K_1+P_3$, a $W_4$ is the graph $K_1+C_4$, a {\em house} is just the complement of $P_5$, a {\em crown} is the graph $K_1+K_{1,3}$, and an $HVN$ is a $K_4$ together with one more vertex which is adjacent to exactly two vertices of $K_4$. (See Figure \ref{fig-1} for these configurations.)

In \cite{PA22}, Prashant {\em  et al.} proved that if $G$ is ($P_3\cup P_2$, diamond)-free, then $\chi(G)\le4$ when $\omega(G)=2$, $\chi(G)\le6$ when $\omega(G)=3$, $\chi(G)\le4$ when $\omega(G)=4$, and $G$ is perfect when $\omega(G)\ge5$, and they also proved \cite{PA22} that $\chi(G)\le\omega(G)+1$ if $G$ is a $(P_3\cup P_2, HVN)$-free graph with $\omega(G)\ge4$.  As a superclass of ($P_3\cup P_2$, diamond)-free graphs, Cameron {\em  et al.} \cite{CHM21} proved that $\chi(G)\le\omega(G)+3$ if $G$ is ($P_6$, diamond)-free, this bound is optimal. In \cite{WZ22}, Wang and Zhang proved that if $G$ is a $(P_3\cup P_2, K_3)$-free
graph, then $\chi(G)\le3$ unless $G$ is one of eight graphs with $\Delta(G)=5$ and $\chi(G)=4$, and they proved \cite{WZ22} $\chi(G)\le3\omega(G)$ if $G$ is $(P_3\cup P_2, W_4)$-free.  Recently, Wu and Xu \cite{WX23} proved that $\chi(G)\le\frac{1}{2}\omega^2(G)+\frac{3}{2}\omega(G)+1$ if $G$ is ($P_3\cup P_2$, crown)-free. Char and Karthick \cite{CK22} proved that $\chi(G)\le$ max $\{\omega(G)+3, \lfloor \frac{3\omega(G)}{2} \rfloor-1\}$ if $G$ is a ($P_3\cup P_2$, $\bar{P_3\cup P_2}$)-free graph with $\omega(G)\ge3$, this bound is optimal.

In this paper, the main result is that

\begin{theorem}\label{P3}
	$\chi(G)\le2\omega(G)$ if $G$ is $($$P_3\cup P_2$, house$)$-free.
\end{theorem}

 Let $G$ be a graph on $n$-vertices $\{v_1, v_2, \dots , v_n\}$ and let $H_1, H_2, \dots , H_n$ be $n$ vertex-disjoint graphs. An expansion
$G(H_1, H_2, \dots , H_n)$ of $G$ is a graph obtained from $G$ by (i) replacing each $v_i$ of $G$ by $H_i$, $i = 1, 2, . . . , n$ and
(ii) by joining every vertex in $H_i$ with every vertex in $H_j$, whenever $v_i$ and $v_j$ are adjacent in $G$. In addition, $G\cong K_n$ and $H_i\cong H$, denote by $K_n(H)=G(H,H,\dots, H)$.

Let $H$ be the Mycielski-Gr\"{o}stzsch graph (see Figure 2). Then $\omega(H)=2$ and $\chi(H)=4$. It is clear that $K_k(H)$ is ($P_3\cup P_2$, house)-free, and $\chi(K_k(H))=2\omega(K_k(H))=2\cdot 2k=4k$. Let $H'$ be the complement of Schl\"{a}fli graph (see https://houseofgraphs.org/graphs/19273). Then $\omega(H')=3$ and $\chi(H)=6$. Obviously, $K_{k-1}(H)+H'$ is ($P_3\cup P_2$, house)-free. Hence $\omega(K_{k-1}(H)+H')=2k+1$ and $\chi(K_{k-1}(H)+H')=4k+2$. This implies that our bound is optimal when $\omega(G)\ge2$.

\begin{figure}[htbp]\label{fig-2}
	\begin{center}
		\includegraphics[width=4cm]{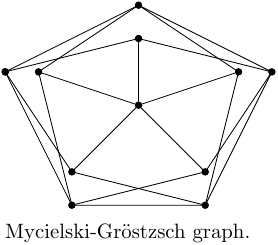}
	\end{center}
	\vskip -15pt
	\caption{Mycielski-Gr\"{o}stzsch graph.}
\end{figure}

\section{Preliminary and notations}

Let $G$ be a graph, $v\in V(G)$, and let $X$ and $Y$ be two subsets of $V(G)$. We say that $v$ is {\em complete} to $X$ if $v$ is adjacent to all vertices of $X$, and say that $v$ is {\em anticomplete} to $X$ if $v$ is not adjacent to any vertex of $X$. We say that $X$ is complete (resp. anticomplete) to $Y$ if each vertex of $X$ is complete (resp. anticomplete) to $Y$.  Particularly, we say that $X$ is {\em almost complete}  to $Y$ if at most one vertex of $X$ is not complete to $Y$. For $u, v\in V(G)$, we simply write $u\sim v$ if $uv\in E(G)$, and write $u\not\sim v$ if $uv\not\in E(G)$.

A {\em hole} of $G$ is an induced cycle of length at least 4, and a {\em $k$-hole} is a hole of length $k$. A $k$-hole is called an {\em odd hole} if $k$ is odd, and is called an {\em even hole} otherwise. An {\em antihole} is the complement of some hole. An odd (resp. even) antihole is defined analogously. The famous {\em Strong Perfect Graph Theorem} states that

\begin{theorem}\label{Perfect}\cite{CRSR06}
	A graph is perfect if and only if it induces neither an odd hole nor an odd antihole.
\end{theorem}

From now on, we may always assume that $G$ is a ($P_3\cup P_2$, house)-free graph such that $\chi(G')\le 2\omega(G')$ for every induced subgraph $G'$ of $G$ different from $G$, and $\chi(G)>2\omega(G)$. Then the next lemma follows easily.

\begin{lemma}\label{nonadjcent}
	Let $u$ and $v$ be two nonadajcent vertices in $G$. Then $N(u)\not\subseteq N(v)$ and $N(v)\not\subseteq N(u)$.
\end{lemma}
\pf Suppose to its contrary that $N(u)\subseteq N(v)$ by symmetry. By assumption, $\chi(G-u)\le 2\omega(G-u) $. Since we can color $u$ by the color of $v$, it follows that $\chi(G)\le 2\omega(G-u) \le 2\omega(G)$, a contradiction.
\qed

\medskip

Let $H$ be an induced subgraph of $G$ such that $\chi(G[H])\le\omega(G[H])$. Let $I=\{v\in V(G) | v$ is anticomplete to $V(H)$$\}$, $R=\{v\in V(G) | v$ is complete to $V(H)$$\}$, and $T=V(G)\setminus (V(H)\cup I\cup R)$. If $\chi(G[I])\le \chi(G[H])$ and $\chi(G[T])\le \chi(G[H])$, then we say $H$ a {\em good} subgraph of $G$ and $G$ has a good subgraph $H$. By a simple induction, we can show that

\begin{lemma}\label{good}
	$G$ has no good subgraphs.
\end{lemma}
\pf Suppose not. Let $H$ be a good subgraph of $G$ with $\omega(H)=\omega_0$, $I=\{v\in V(G) | v$ is anticomplete to $V(H)$$\}$, $R=\{v\in V(G) | v$ is complete to $V(H)$$\}$, and $T=V(G)\setminus (V(H)\cup I\cup R)$. Then $\chi(G)\le \chi(G[V(H)\cup I])+\chi(G[T])+\chi(G[R])\le 2\omega_0+2\omega(G[R])\le2\omega_0+2(\omega(G)-\omega_0)=2\omega(G)$, a contradiction.
\qed

\medskip

Since $P_5$-free graphs and $(P_3\cup P_2)$-free graphs are both superclasses of $2K_2$-free graphs. By the following theorem, we may also always assume that $G$ induces an $2K_2$.

\begin{theorem}\cite{FGM1995}\label{P5}
	If $G$ is a $($$P_5$, house$)$-free  graph, then $\chi(G)\le \lfloor \frac{3\omega}{2} \rfloor$.
\end{theorem}

Let $P=G[\{v_1,v_2,v_3,v_4\}]$ be an induced $2K_2$ such that $v_1\sim v_2$ and $v_3\sim v_4$.  For a subset $S\subseteq\{1,2,3,4\}$, we define $N_S=\{v|v\in N(P),$ and $v_i\sim v$ if and only if $i\in S\}$. Let $N=N(\{v_1,v_2,v_3,v_4\})$ and $A=M(\{v_1,v_2,v_3,v_4\})$. Since $G$ is $(P_3\cup P_2)$-free, we have that $N_S=\emptyset$ if $|S|=1$. Therefore, $N=N_{\{1,2\}}\cup N_{\{1,3\}}\cup N_{\{1,4\}}\cup N_{\{2,3\}}\cup N_{\{2,4\}}\cup N_{\{3,4\}}\cup N_{\{1,2,3\}}\cup N_{\{1,2,4\}}\cup N_{\{1,3,4\}}\cup N_{\{2,3,4\}}\cup N_{\{1,2,3,4\}}$.

\begin{lemma}\label{c0}
	Let $\{i,i'\}=\{1,2\}, \{j,j'\}=\{3,4\}$. Then $N_{\{i,j\}}$ is anticomplete to $N_{\{i,j'\}}\cup N_{\{i',j\}}
	\cup N_{\{i,i',j'\}}\cup N_{\{i',j,j'\}}$, and is complete to $N_{\{i,i',j\}}\cup N_{\{i,j,j'\}}\cup N_{\{1,2,3,4\}}$. Furthermore, let $\{s,s'\}\in \{\{1,2\},\{3,4\}\}$ and $\{t,t'\}=\{1,2,3,4\}\setminus \{s,s'\}$. Then $N_{\{s,s',t\}}$ is anticomplete to $N_{\{s,s',t'\}}$, and is complete to $N_{\{s,t,t'\}}\cup N_{\{s',t,t'\}}\cup N_{\{1,2,3,4\}}$.
\end{lemma}
\pf  By symmetry, we may assume that $i=1,j=3$, which means $i'=2,j'=4$. Suppose there exist two adjacent vertices $u_1,u_2$ such that $u_1\in N_{\{1,3\}}$ and $u_2\in N_{\{2,3\}}$. Then $\{v_1,v_2,v_3,u_1,u_2\}$ induces a house, a contradiction. So, $N_{\{1,3\}}$ is anticomplete to $N_{\{2,3\}}$. Similarly, $N_{\{1,3\}}$ is anticomplete to $N_{\{1,4\}}\cup N_{\{1,2,4\}}\cup N_{\{2,3,4\}}$.
Suppose there exist two nonadjacent vertices $u_1',u_2'$ such that $u_1'\in N_{\{1,3\}}$ and $u_2'\in N_{\{1,2,3\}}$. Then $\{v_1,v_2,v_3,u_1',u_2'\}$ induces a house, a contradiction. So, $N_{\{1,3\}}$ is complete to $N_{\{1,2,3\}}$. Similarly, $N_{\{1,3\}}$ is complete to $N_{\{1,3,4\}}\cup N_{\{1,2,3,4\}}$.

By the same arguments as above, we can easily deduce the latter conclusion.
\qed

\medskip

We will give some structural properties of $G$ in Section 3, and prove Theorem \ref{P3} in Section 4.

\section{Structural properties of $G$}

In this section, we  will give some structural properties of $G$.

\medskip

\textbf{(F1)} Let $|S|=3$. Then $G[N_S]$ is $P_3$-free.

\medskip

\pf By symmetry, let $ u_1$-$u_2$-$u_3$ be an induced $P_3$ in $N_{\{1,3,4\}}$. By Lemma \ref{nonadjcent}, there exists a vertex $v\in N(v_2)$ such that $v\not\sim u_3$.

Suppose that $v\not\sim v_1$. Then  $v\sim v_3$ or $v\sim v_4$ as otherwise $\{v_1,v_2,v,v_3,v_4\}$ induces a $P_3\cup P_2$, a contradiction. By Lemma \ref{c0}, $v$ is anticomplete to $\{u_1,u_2,u_3\}$, which implies that $\{u_1,u_2,u_3,v,v_2\}$ induces a $P_3\cup P_2$, a contradiction.

Thus $v\sim v_1$. If $v\sim v_3$ , then $\{u_3,v_3,v, v_1,v_2\}$ induces a house, a contradiction. By symmetry, $v\not\sim v_4$. Then $v\sim u_1$ as otherwise $\{u_1,v_4,u_3,v_2,v\}]$ induces a $P_3\cup P_2$, a contradiction. But now, $\{v_1,v,u_1,v_3,u_3\}$ induces a house, a contradiction.\qed

\medskip

\textbf{(F2)} Let $|S|=2$. Then $G[N_S]$ is $P_3$-free.

\medskip

\pf Since $G$ is $P_3\cup P_2$-free, we have that $G[N_{\{1,2\}}]$ and $G[N_{\{3,4\}}]$ are both $P_3$-free. By symmetry, we may assume that there exists an induced $P_3$ in $G[N_{\{1,3\}}]$, say $u_1$-$u_2$-$u_3$. Since $u_1\not\sim v_2$, there exists a vertex $v\in N(v_2)$ such that $v\not\sim u_1$. By Lemma \ref{c0}, $v\in N_{\{1,2\}}\cup N_{\{2,3\}}\cup N_{\{2,4\}}\cup N_{\{1,2,4\}}\cup N_{\{2,3,4\}}$.

Note that $v$ has at least two neighbors in $V(P)$. If $v\in N_{\{1,2\}}$, then $\{u_1,v_3,v_4,v_2,v\}$ induces a $ P_3\cup P_2$, a contradiction. So, $v\not\in N_{\{1,2\}}$. If $v\in N_{\{2,3\}}\cup N_{\{1,2,4\}}\cup N_{\{2,3,4\}}$, by Lemma \ref{c0}, $v$ is anticomplete to $N_{\{2,3\}}\cup N_{\{1,2,4\}}\cup N_{\{2,3,4\}}$, then $\{u_1,u_2,u_3,v,v_2\}$ induces a $P_3\cup P_2$, a contradiction. So, $v\in N_{\{2,4\}}$. To forbid an induced $P_3\cup P_2$ on $\{u_1,u_2,u_3,v,v_2\}$, we have that $v\sim u_2$ or $v\sim u_3$.

If $v\sim u_2$, then $\{u_1,u_2,v_3,v_4,v\}$ induces a house, a contradiction. So, $v\not\sim u_2$. If $v\sim u_3$, then $\{u_2,u_3,v_3,v_4,v\}$ induces a house, a contradiction.\qed

\medskip

\textbf{(F3)} $N_{\{1,2\}}$ and $N_{\{3,4\}}$ are both cliques.

\medskip

\pf Suppose there exist two nonadjacent vertices $u$ and $v$ in $N_{\{1,2\}}$. Then $\{u,v_1,v,v_3,v_4\}$ induces a $P_3\cup P_2$, a contradiction. So, $N_{\{1,2\}}$ is a clique, and $N_{\{3,4\}}$ is a clique by symmetry.

\medskip

\textbf{(F4)} $A$ is anticomplete to $N_{\{1,2\}}\cup N_{\{3,4\}}$.

\medskip

\pf Suppose not. By symmetry, there exist two adjacent vertices $x$ and $y$ such that $x\in N_{\{1,2\}}$ and $y\in A$. Then $\{v_1,x,y,v_3,v_4\}$ induces a $P_3\cup P_2$, a contradiction. So, $A$ is anticomplete to $N_{\{1,2\}}$, and is anticomplete to $N_{\{3,4\}}$ by symmetry.

\medskip

\textbf{(F5)} Let $|S|=3$. If $G[N_S]$ induces a $P_2\cup P_1$, then $G[N_S]$ is a good subgraph of $G$.

\medskip

\pf By symmetry, let $S=\{1,3,4\}$. By \textbf{(F1)}, we have that $G[N_{\{1,3,4\}}]$ is  an union of cliques. Since $G[N_{\{1,3,4\}}]$ induces a $P_2\cup P_1$, we have that $G[N_{\{1,3,4\}}]$ has two cliques, say $B$ and $C$, and $|C|\ge 2$ by symmetry. Let $u_1,u_2\in C$ and $u_3\in B$. By Lemma \ref{c0}, we have that $N_{\{1,3,4\}}$ is complete to $N_{\{1,3\}}\cup N_{\{1,4\}}\cup N_{\{1,2,3\}}\cup N_{\{1,2,4\}}\cup N_{\{1,2,3,4\}}$, and is anticomplete to $N_{\{2,3\}}\cup N_{\{2,4\}}\cup N_{\{2,3,4\}}$.

First, we prove that $N_{\{2,3\}}=N_{\{2,4\}}= \emptyset$. If $N_{\{2,3\}}\ne \emptyset$, let $v\in N_{\{2,3\}}$, then $v$ is anticomplete to $\{u_1,u_2,u_3\}$ by Lemma \ref{c0}, which implies that $\{v_2,v\}\cup \{u_3,v_4,u_1\}$ induces a $ P_3\cup P_2$, a contradiction. Therefore, $N_{\{2,3\}}=\emptyset$. By symmetry, $N_{\{2,4\}}= \emptyset$.

Second, we prove that $|N_{\{2,3,4\}}|\le 1$. By Lemma \ref{c0}, we have that $N_{\{2,3,4\}}$ is anticomplete to $N_{\{1,3,4\}}$. Suppose that $v,v'$ are two distinct vertices in $N_{\{2,3,4\}}$. Then
$\{v,v_2,v',u_1,u_2\}$ induces a $P_3\cup P_2$ if $v\not\sim v'$, and $\{v_1,u_1,v_3,v,v'\}$ induces a $P_3\cup P_2$ if $v\sim v'$, both are contradictions. Thus $|N_{\{2,3,4\}}|\le 1$.

Third, we prove that $N_{\{1,2\}}$ is complete to $B\cup C$. For each $v\in N_{\{1,2\}}$, we have that $v\sim u_1$ or $v\sim v_3$ as otherwise $\{u_1,v_3,u_3,v_2,v\}$ induces a $P_3\cup P_2$, a contradiction. If $v\not\sim u_1$ or $v\not\sim u_3$, then $\{u_3,u_1,v_3,v,v_1\}$ induces a house, a contradiction. Therefore, $v\sim u_1$ and $v\sim u_3$. If $v\not\sim u_2$, then $\{u_3,u_2,v_3,v,v_1\}$ induces a house, a contradiction. So, $N_{\{1,2\}}$ is complete to $B\cup C$.

 Next, we prove that $N_{\{3,4\}}$ is complete to $C$. If there are two nonadjacent vertices $v,v'$ in $N_{\{3,4\}}$, then $\{v,v_3,v',v_1,v_2\}$ induces a $P_3\cup P_2$, a contradiction. It follows that $N_{\{3,4\}}$ is a clique. By symmetry, we may assume that $v\not\sim u_1$. To forbid an induced house on $\{v,v_3,u_3,v_1,u_1\}$, we have that $v\not\sim u_3$. By Lemma \ref{nonadjcent}, there exists a vertex $v'\in N(v)$ such that $v'\not\sim u_3$, which implies that $v\not\sim u_1$ as otherwise $\{v,u_3,v_1 ,u_1,v_3\}$ induces a house, a contradiction. To forbid an induced  $P_3\cup P_2$ on $\{u_1,v_1,v_2,v,v'\}$,  we have that $v'\sim v_1$ or $v'\sim v_2$. So $v'\in N_{\{2,3,4\}}$. But now, $\{u_1,v_1,u_3,v,v'\}$ induces a $P_3\cup P_2$, a contradiction. Hence $N_{\{3,4\}}$ is complete to $C$.

 Now, let $A':=\{v\ | v\in A$ and $v$ is not complete to $C\}$, we prove that $A'$ is an independent set.
 Suppose there exist two adjacent vertices $u,v$ in $G[A]$ such that $u$ and $v$ are both not complete to $C$. For any $i\in\{1,2,3\}$, if $u\not\sim u_i$ and $v\not\sim u_i$, then $\{u_i,v_1,v_2,u,v\}$ induces a $P_3\cup P_2$, a contradiction. So, $|N(u_i)\cap \{u,v\}|\ge1$ for any $i\in\{1,2,3\}$. Thus, we may suppose that $u\not\sim u_1$ and $v\not\sim u_2$. Then $u\sim u_2$ and $v\sim u_1$. But now, $\{v_1,u_1,u_2,u,v\}$ induces a house, a contradiction. So, $A'$ is an independent set.

 Let $\omega(N_{\{1,3,4\}})=\omega_0$, and $T=\{v\in V(G) | v$ is not complete to $N_{\{1,3,4\}}$$\}$. Then $T=A'\cup \{v_2\}\cup N_{\{2,3,4\}}$. Since $v_2$ is anticomplete to $A'$, and $N_{\{2,3,4\}}$ has at most one vertex, we have that $\chi(G[T])\le2\le \omega_0$, which implies that $G[N_{\{1,3,4\}}]$ is a good subgraph of $G$.
\qed

\medskip

\textbf{(F6)} Let $|S|=3$ and $A\neq \emptyset$. If $N_S$s are all cliques, then there exists an $S$ such that $G[N_S]$ is a good subgraph of $G$.

\medskip

\pf By \textbf{(F5)}, we may assume that all the $N_Ss$ contain no $P_2\cup P_1$ when $|S|=3$, that is, $N_S$s are all cliques. Choose one with maximum cardinality from the $N_Ss$. Without loss of generality, we consider $S=\{1,3,4\}$, that is to say, $N_{\{1,3,4\}}$ is a clique with $\omega_0=|N_S|\geq2$. By Lemma \ref{c0}, we have that $N_{\{1,3,4\}}$ is complete to $N_{\{1,3\}}\cup N_{\{1,4\}}\cup N_{\{1,2,3\}}\cup N_{\{1,2,4\}}\cup N_{\{1,2,3,4\}}$, and is anticomplete to $N_{\{2,3\}}\cup N_{\{2,4\}}\cup N_{\{2,3,4\}}$.

First, we prove that $N_{\{3,4\}}$ is almost complete to $N_{\{1,3,4\}}$. Suppose that there exist two vertices in $N_{\{3,4\}}$ are not complete to $N_{\{1,3,4\}}$, say $x,y$. If $x\not\sim y$, then $\{x,v_4,y,v_1,v_2\}$ induces a $P_3\cup P_2$, a contradiction. So, $x\sim y$. For each $u'\in N_{\{1,3,4\}}$, we see that $u'\sim x$ or $u'\sim y$ as otherwise $\{u',v_1,v_2,x,y\}$ induces a $P_3\cup P_2$, a contradiction. Therefore, we may assume that there exist two distinct vertices $u_1,u_2$ in $N_{\{1,3,4\}}$ such that $u_1\not\sim x$ and $u_2\not\sim y$. Then $u_1\sim y$ and $u_2\sim x$. But now, $\{v_1,u_1,x,y,u_2\}$ induces a house, a contradiction. It implies that at most one vertex in $N_{\{3,4\}}$ is not complete to $N_{\{1,3,4\}}$.

Second, we prove that $N_{\{1,3,4\}}$ is complete to $A$. Let $x\in A$. Assume that $x\not\sim u_1$. Then there exists a vertex $y\in N(x)$ such that $y\not\sim u_1$ by Lemma \ref{nonadjcent},. Note that $y\in N(P)$, for otherwise, $\{v_1,u_1,v_3,x,y\}$ induces a $P_3\cup P_2$, a contradiction. Hence, $y\in N_{\{1,2\}}\cup N_{\{3,4\}}\cup N_{\{2,3\}}\cup N_{\{2,4\}}\cup N_{\{2,3,4\}}$. If $y\in N_{\{1,2\}}$, then $\{x,y, v_2,v_3,v_4\}$ induces a $P_3\cup P_2$, a contradiction. If $y\in N_{\{3,4\}}$, then $\{x,y, v_4,v_1,v_2\}$ induces a $P_3\cup P_2$, a contradiction. Since $N_{\{1,3,4\}}$ is anticomplete to $N_{\{2,3\}}\cup N_{\{2,4\}}\cup N_{\{2,3,4\}}$, we have that $y$ is anticomplete to $u_1,u_2$ if $y\in N_{\{2,3\}}\cup N_{\{2,4\}}\cup N_{\{2,3,4\}}$. But now, $\{x,y,v_2,u_1,u_2\}$ induces a $P_3\cup P_2$, a contradiction. So, $N_{\{1,3,4\}}$ is complete to $A$.

Third, we prove that $N_{\{1,2\}}$ is complete to $N_{\{1,3,4\}}$. Suppose not. Let $v\in N_{\{1,2\}}$ such that $v\not\sim u_1 $. Since $N_{\{1,3,4\}}$ is complete to $A$ and $A\ne\emptyset$, let $x\in A$, we have that $x\sim u_1$. Then $x\not\sim v$, for otherwise, $[\{v_1,v_2,v,x,u_1\}$ induces a house, a contradiction. But now, $\{x,u_1,v_3\,v_2,v\}$ induces a $P_3\cup P_2$, a contradiction. So, $N_{\{1,2\}}$ is complete to $N_{\{1,3,4\}}$.

Next, we prove that $|N_{\{2,3\}}|\leq 1$ and $|N_{\{2,4\}}|\leq 1$. Suppose not. We may by symmetry assume that there are two vertices in $N_{\{2,3\}}$, say $z_1, z_2$. Then $\{v_1,u_1,v_4,z_1, z_2\}$ induces a $P_3\cup P_2$ if $z_1\sim z_2$ and $\{z_1,v_2,z_2,u_1,u_2\}$ induces a $P_3\cup P_2$ if $z_1\not\sim z_2$, both are contradictions. So, $|N_{\{2,3\}}|\leq 1$, and $|N_{\{2,4\}}|\leq 1$ by symmetry.

Suppose there exists a vertex in $N_{\{3,4\}}$ which is not complete to $N_{\{1,3,4\}}$, say $z$. Note that $\chi(G[\{z,v_2\}])]\le1$ as $\{z,v_2\}$ is an independent set. Since $N_{\{2,3\}}\cup N_{\{2,4\}}$ is an independent set by Lemma \ref{c0} and $|N_{\{2,3\}}|\leq 1, |N_{\{2,4\}}|\leq 1$, we have that $\chi(G[N_{\{2,3\}}\cup N_{\{2,4\}}])\le1$.  Furthermore, $N_{\{2,3,4\}}$ is anticomplete to $N_{\{1,3,4\}}$, and $V(G)\setminus(V(H)\cup \{z,v_2\}\cup N_{\{2,3\}}\cup N_{\{2,4\}})$ is complete to $N_{\{1,3,4\}}$. So, $G[N_{\{1,3,4\}}]$ is a good subgraph of $G$ as $\omega(G[N_{\{2,3,4\}}])\le\omega_0$ and $\chi(G[\{z,v_2\}\cup N_{\{2,3\}}\cup N_{\{2,4\}}])\le2$.
\qed

\medskip

\textbf{(F7)} Let $i\in\{1,2\}$ and $j\in\{3,4\}$. If $N_{\{i,j\}}$ induces a $P_2\cup P_1$, then $G[N_{\{i,j\}}]$ is a good subgraph of $G$.

\medskip

\pf By symmetry, let $i=1$ and $j=3$. By \textbf{(F2)}, we have that $N_{\{1,3\}}$ is $P_3$-free, which implies that it is an union of cliques. We choose a maximal clique in $N_{\{1,3\}}$, say $C$. Since $N_{\{i,j\}}$ induces a $P_2\cup P_1$, we have that $|C|=\omega_0\ge2$ and there exists a clique $B$ other than $C$. Let $u_1,u_2\in C$ and $u_3\in B$. By Lemma \ref{c0}, we have that $N_{\{1,3\}}$ is complete to $N_{\{1,2,3\}}\cup N_{\{1,3,4\}}\cup N_{\{1,2,3,4\}}$, and is anticomplete to $N_{\{2,3\}}\cup N_{\{1,4\}}\cup N_{\{1,2,4\}}\cup N_{\{2,3,4\}}$.

First, we prove that $N_{\{2,4\}}$ is complete to $C$. Let $v\in N_{\{2,4\}}$. If $|C\setminus N(v)|\ge2$, then $\{v_2,v,v_3,x_1,x_2\}$ induces a $P_3\cup P_2$, where $\{x_1,x_2\}\subseteq C\setminus N(v)$, a contradiction. So, $|C\setminus N(v)|\le1$. If $|C\setminus N(v)|=1$, then there exist two vertices $u_1,u_2$ in $C$ such that $u_1\sim v$ and $u_2\not\sim v$ as $|C|\ge2$. But now, $\{v_1,v_2,v,u_1,u_2\}$ induces a house, a contradiction. Therefore, $N_{\{2,4\}}$ is complete to $C$.

Second, we prove that $|N_{\{1,4\}}|\le1$ and $|N_{\{2,3\}}|\le1$. Suppose to its contrary that we may by symmetry assume that $|N_{\{1,4\}}|\ge2$. Let $x_1,x_2\in N_{\{1,4\}}$. Suppose $x_1\not\sim x_2$. Since $N_{\{1,4\}}$ is anticomplete to $N_{\{1,3\}}$, we have that $\{x_1,x_2\}$ is anticomplete to $\{u_1,u_2\}$. But now,  $\{x_1,x_2,v_4,u_1,u_2\}$ induces a $P_3\cup P_2$ if $x_1\not\sim x_2$ and $\{x_1,x_2,v_3,u_1,u_3\}$ induces a $P_3\cup P_2$ if $x_1\sim x_2$, both are contradictions.  So,
$|N_{\{1,4\}}|\le1$, and $|N_{\{2,3\}}|\le1$ by symmetry.

In reality, $N_{\{1,4\}}\cup N_{\{2,3\}}$ is an independent set. Suppose not, let $\{x_1'\}=N_{\{1,4\}}$ and $\{x_2'\}=N_{\{2,3\}}$ such that $x_1'\sim x_2'$. Since $N_{\{1,3\}}$ is anticomplete to $N_{\{2,3\}}\cup N_{\{1,4\}}$, we have that $\{u_1,u_2\}$ is anticomplete to $\{x_1',x_2'\}$. But now, $\{x_1',x_2'v_2,u_1,u_2\}$ induces a $P_3\cup P_2$, a contradiction.

Third, we prove that $N_{\{1,2\}}\cup N_{\{3,4\}}$ is complete to $B\cup C$. For each $v\in N_{\{1,2\}}$, we have that $v\sim u_1$ or $v\sim u_3$ as otherwise $\{u_1,v_3,u_3,v_2,v\}$ induces a $P_3\cup P_2$, a contradiction. If $v\not\sim u_1$ or $v\not\sim u_3$, then $\{u_3,u_1,v_3,v,v_1\}$ induces a house, a contradiction. Therefore, $v\sim u_1$ and $v\sim u_3$. If $v\not\sim u_2$, then $\{u_3,u_2,v_3,v,v_1\}$ induces a house, a contradiction. So, $N_{\{1,2\}}$ is complete to $B\cup C$.

Next, we prove that $N_{\{1,2,4\}}= N_{\{2,3,4\}}=\emptyset$. By symmetry, we suppose that there exists a vertex $x_1$ in $N_{\{1,2,4\}}$. If $x_1$ is anticomplete to $C$, then $\{x_1,v_2,v_4,u_1,u_2\}$ induces a $P_3\cup P_2$, a contradiction. If $x_1$ is not anticomplete to $C$, without loss of generality, $x_1\not\sim u_1$, then $\{x_1,v_1,u_1,v_3,v_4\}$ induces a house, a contradiction. So, $N_{\{1,2,4\}}=\emptyset$, and $N_{\{2,3,4\}}=\emptyset$ by symmetry.

Let $A'=\{v\ | v\in A$ and $v$ is not complete to $C\}$. By the same arguments in \textbf{(F5)}, we have that $A'$ is an independent set. Furthermore, $N_{\{1,4\}}\cup N_{\{2,3\}}$ is an independent set.
Therefore, $G[N_{\{1,3\}}]$ is a good subgraph of $G$ as $\chi(G[A'\cup N_{\{1,4\}}\cup N_{\{2,3\}}])\le2$ and $V(G)\setminus(N_{\{1,3\}}\cup A'\cup N_{\{1,4\}}\cup N_{\{2,3\}})$ is complete to $N_{\{1,3\}}$.\qed

\medskip

\textbf{(F8)} Let $i\in\{1,2\}$ and $j\in\{3,4\}$. If some $N_{\{i,j\}}$ is a clique with $\omega(G[N_{\{i,j\}}])=\omega_0\ge2$, then $G[N_{\{i,j\}}]$ is a good subgraph of $G$.

\medskip

\pf  Suppose that $N_{\{1,3\}}$ is a clique with order at least two. Let $u_1,u_2\in N_{\{1,3\}}$. By Lemma \ref{c0}, we have that $N_{\{1,3\}}$ is complete to $N_{\{1,2,3\}}\cup N_{\{1,3,4\}}\cup N_{\{1,2,3,4\}}$, and is anticomplete to $N_{\{2,3\}}\cup N_{\{1,4\}}\cup N_{\{1,2,4\}}\cup N_{\{2,3,4\}}$.

First, we prove that $N_{\{1,3\}}$ is complete to $N_{\{2,4\}}$. Suppose to its contrary that  $x\in N_{\{2,4\}}$ such that $x\not\sim u_1$.  To forbid an induced $P_3\cup P_2$ on $\{v_2,x,v_4,u_1,u_2\}$, we have that $x\sim u_2$. But now, $\{v_1,v_2,x,u_1,u_2\}$ induces a house, a contradiction. So, $N_{\{1,3\}}$ is complete to $N_{\{2,4\}}$.

Second, we prove that  $N_{\{1,4\}}$  is anticomplete to $N_{\{2,3\}}$, and either they are both singletons or one of them is empty. Let $x\in N_{\{1,4\}}$, and $y\in N_{\{2,3\}}$. Then $x\not\sim y$ as otherwise $\{x,y,v_2,u_1,u_2\}$ induces a $P_3\cup P_2$, a contradiction. If $|N_{\{1,4\}}|\geq2$, let $x'$ be another vertex in $N_{\{1,4\}}$, then $\{v_2,y,v_3,x,x'\}$ induces a $P_3\cup P_2$, a contradiction. So, $N_{\{1,4\}}$  is anticomplete to $N_{\{2,3\}}$, and either they are both singletons or one of them is empty.

Third, we prove $ N_{\{2,3,4\}}\cup N_{\{1,2,4\}}=\emptyset$. Suppose there exists a vertex $v\in N_{\{2,3,4\}}$. Then $\{v_2,v,v_4,u_1,u_2\}$ induces a $P_3\cup P_2$, a contradiction. So, $ N_{\{2,3,4\}}=\emptyset$, and $N_{\{1,2,4\}}=\emptyset$ by symmetry.

Next, we prove that $N_{\{1,3\}}$ is complete to $N_{\{1,2\}}\cup N_{\{3,4\}}\cup A$. Suppose that there exists a vertex $v\in N_{\{1,2\}}$ such that $v\not\sim u_1$. Then $\{u_1,v_3,v_4,v,v_2\}$ induces a $P_3\cup P_2$, a contradiction. Thus $N_{\{1,3\}}$ is complete to $N_{\{1,2\}}$, and also complete to $N_{\{3,4\}}$ by symmetry. Suppose there exists a vertex $x\in A$ such that $x\not\sim u_1$. Then $x$ has a neighbor $y$ such that $y\not\sim u_1$. Naturally, $y\in N(P)$. Otherwise, $\{v_1,u_1,v_3,x,y\}$ induces a $P_3\cup P_2$, a contradiction. Hence $y\in N_{\{2,3\}}\cup N_{\{1,4\}} $. If $y\in N_{\{2,3\}}$, then $x\not\sim u_2$ as otherwise $\{u_1,u_2,x,y,v_3\}$ induces a house, a contradiction. But then $\{x,y,v_2,u_1,u_2\}$ induces a $P_3\cup P_2$, a contradiction. Therefore, $N_{\{1,3\}}$ is complete to $A$.

Note that $\chi(G[N_{\{2,3\}}\cup N_{\{1,4\}}\cup N_{\{1,2,4\}}\cup N_{\{2,3,4\}}])\le\omega_0$ as $ N_{\{2,3,4\}}\cup N_{\{1,2,4\}}=\emptyset$, $N_{\{1,4\}}$  is anticomplete to $N_{\{2,3\}}$, and $\omega(G[N_{\{1,4\}}])\le \omega_0, \omega(G[N_{\{2,3\}}])\le \omega_0$.  Furthermore, $\chi(\{v_2,v_4\})\le1$. The rest vertices are complete to $N_{\{1,3\}}$. Therefore,  $G[N_{\{1,3\}}]$ is a good subgraph of $G$.\qed

\section{Proof of Theorem \ref{P3}}

In this section, we will complete the proof of Theorem \ref{P3}. Let $B$ be the union of the sets $N_S$ for $|S|=3$ and the sets $N_{\{i,j\}}$ for $i\in\{1,2\}$ and $j\in\{3,4\}$. That is $B=(N_{\{2,3\}}\cup N_{\{2,4\}}\cup N_{\{1,3,4\}})\cup (N_{\{1,3\}}\cup N_{\{1,4\}}\cup N_{\{2,3,4\}})\cup (N_{\{1,2,3\}}\cup N_{\{1,2,4\}}$). By Lemma \ref{c0}, $B$ is complete to $N_{\{1,2,3,4\}}$.

We firstly need the following two lemmas.

\begin{lemma}\label{B}
	$\chi(G[B])\le3$ if $A\ne\emptyset$, and $G[B]$ is perfect if $A=\emptyset$.
\end{lemma}

\pf Suppose that $A\ne\emptyset$. By \textbf{(F5)}, \textbf{(F6)}, \textbf{(F7)}, \textbf{(F8)} and Lemma \ref{c0}, we have that $N_{\{2,3\}}\cup N_{\{2,4\}}\cup N_{\{1,3,4\}}$, $N_{\{1,3\}}\cup N_{\{1,4\}}\cup N_{\{2,3,4\}}$ and $N_{\{1,2,3\}}\cup N_{\{1,2,4\}}$ are all independent, which implies that $\chi(G[B])\le3$.

Suppose that $A=\emptyset$. By \textbf{(F1)}, we have that each $N_S$ is either an independent set or a clique with $|S|=3$. By  \textbf{(F7)} and  \textbf{(F8)}, we have that each $N_{\{i,j\}}$ is an independent set for $i\in\{1,2\}$ and $j\in\{3,4\}$. Consequently, by Lemma \ref{c0}, we may express the structure of $G[B]$ by Figure 3, the solid adjacent sets are complete, nonadjacent sets are anticomplete, dotted adjacent sets are not complete.

If each $N_S$ is not empty with $|S|=3$, then $B=\mathop{\bigcup}\limits_{|S|=3} N_S$ as $G$ is house-free, hence $\chi(G[B])=\omega(G[B])$. If some $N_S$ is empty with $|S|=3$, say $N_{\{2,3,4\}}=\emptyset$, then it is easy to verify that neither $ G[B] $ nor $\overline{G}[B]$ has odd hole. Hence $G[B]$ is perfect. \qed

\begin{figure}[htbp]\label{fig-3}
	\begin{center}
		\includegraphics[width=6cm]{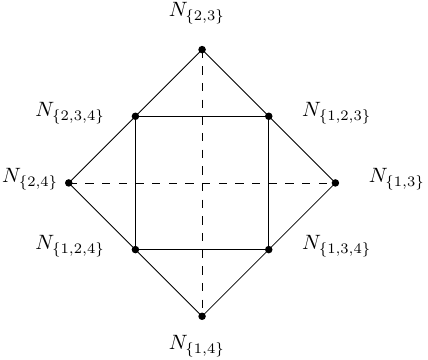}
	\end{center}
	\vskip -15pt
	\caption{Structure of $G[B]$.}
	{The solid adjacent sets are complete, nonadjacent sets are anticomplete, dotted adjacent sets are not complete.}
\end{figure}

\begin{lemma}\label{c7}
	If $N_{\{1,2,3,4\}} \neq \emptyset$, then $N_{\{1,2,3,4\}}$ is a clique.
\end{lemma}

\pf Suppose to its contrary. Let $D=\{v\in N_{\{1,2,3,4\}}| v$ is complete to $ N_{\{1,2,3,4\}}\setminus \{v\} \}$, and $C=N_{\{1,2,3,4\}}\setminus D$.
We can observe that each vertex in $C$ has a nonadjacent vertex in $C$, and $D$ is complete to $C$, and $D$ is a clique. Note that $N_{\{1,2\}}$ and $N_{\{3,4\}}$ are both cliques, and $A$ is anticomplete to $N_{\{1,2\}}\cup N_{\{3,4\}}$, by \textbf{(F3)} and \textbf{(F4)}.

We will prove that $N_{\{1,2\}}$ is almost complete to $C$ and $N_{\{3,4\}}$ is almost complete to $C$, and each vertex of  $N_{\{1,2\}}\cup N_{\{3,4\}}$ is either complete to  $C$ or anticomplete to $C$.

Suppose $x_0,x_1\in N_{\{1,2\}}$ such that $x_0$ and $x_1$ are both not complete to $C$. Suppose $u\in C$ is anticomplete to $\{x_0,x_1\}$. Since $u$ has a nonadjacent vertex in $C$, say $u_1$, we have that $x_0\sim u_1$ or $x_1\sim u_1$ as otherwise $\{u,v_3,u_1,x_0,x_1\}$ induces a $P_3\cup P_2$, a contradiction. But now, $\{u,v_2,u_1,v_3,x_0\}$ induces a house if $x_0\sim u_1$ and $\{u,v_2,u_1,v_3,x_1\}$ induces a house if $x_1\sim u_1$, both are contradictions. So, $x_0$ and $x_1$ cannot be both adjacent to some vertex in $C$.

Since $x_0$ and $x_1$ are both not complete to $C$, let $x_0\not\sim u'$ and $x_1\not\sim u_1'$, it follows that $x_0\sim u_1'$ and $x_1\sim u'$. But now, $\{u,v_2,u_1,v_3,x_0\}$ induces a house if $u'\not\sim u_1'$ and $\{x_0,x_1,u',u_1',v_4\}$ induces a house if $u'\sim u_1'$, both are contradictions. So, $C$ is complete to $N_{\{1,2\}} \setminus \{x_0\}$.

Suppose $x_0$ is either complete to $C$ or anticomplete to $C$. Let $u_2,u_2'\in C$ such that $x_0\not\sim u_2$ and $x_0\sim u_2'$. Then $x$ has a neighbor $y$ such that $y\not\sim u_2$. If $y\in C$, then $\{v_2,v_3,u_2,y,x\}$ induces a house, a contradiction. Hence $y\in N_{\{1,2\}}\cup N_{\{3,4\}}$ as $B$ is complete to $N_{\{1,2,3,4\}}$.  Since $N_{\{1,2\}} \setminus \{x_0\}$ is complete to $C$, we have that $y\in N_{\{3,4\}} $. But now, $\{u_2',v_3,y,x_0,v_2\}$ induces a house, a contradiction. So, $x_0$ is anticomplete to $C$ if $x_0$ is not complete to $C$.

By symmetry, $C$ is complete to $N_{\{3,4\}} \setminus \{y_0\}$ for some $y_0 \in N_{\{3,4\}} $, and $y$ is anticomplete to $C$ if $y_0$ is not complete to $C$.

Now, let $\omega(G[N_{\{1,2,3,4\}}])=\omega_0$ and $\omega(G[C])=\omega_1$, we divide the proof process into three cases depending on $\omega(G[A])\geq 3$ or $1\le\omega(G[A])\leq 2$, or $A=\emptyset$.

\medskip

{\bf Case 1} $\omega(G[A])\geq 3$.

\medskip

Let $Q\subseteq A$ be a clique with $|Q|\geq 3$. We will prove that $Q$ is complete to $N_{\{1,2,3,4\}}$.  We claim that $|Q\setminus N(u)|\le1$. Otherwise, let $a,b\in Q$ such that $a\not\sim u, b\not\sim u$. Then $\{v_1,u,v_3,a,b\}$ induces a $P_3\cup P_2$, a contradiction. So, $|Q\setminus N(u)|\le1$. Suppose there exists a vertex $z$ in $Q$ such that $z\not\sim u$. Let $x,y\in Q$ with $x\ne z, y\ne z$. Then $x\sim u, y\sim u$. By Lemma \ref{nonadjcent}, we have that $z$ has a nonadjacent vertex $z'$ such that $z'\not\sim u$. As $A$ is anticomplete to $N_{\{1,2\}}\cup N_{\{3,4\}}$ and
$C$ is complete to $N_{\{1,2\}}\cup N_{\{3,4\}}\cup B$, it follows that $z'\in C$. Consequently, $z'\sim x$ or $z'\sim y$. By symmetry, $z'\sim x$. Then $\{v_1,u,z',z,x\}$ induces a house, a contradiction. Therefore, $Q$ is complete to $N_{\{1,2,3,4\}}$.

First, we prove that $C$ is complete to $N_{\{1,2\}}\cup N_{\{3,4\}}$. Suppose not. By symmetry, there exist two nonadjacent vertices $x$ and $u$ such that $x\in N_{\{1,2\}}$ and $u\in C$. Then $x$ has a neighbor $y$ such that $y\not\sim u$. If $y\in C$, then $\{v_2,v_3,u,y,x\}$ induces a house, a contradiction. Hence $y\in N_{\{1,2\}}\cup N_{\{3,4\}}$ as $B$ is complete to $N_{\{1,2,3,4\}}$. We first consider $y\in N_{\{1,2\}}$. Let $u'\in C$ and $u\not\sim u'$. Then $u'\not\sim x$. If $u'\sim x$, then $\{v_1,x,u',v_3,u\}$ induces a house, a contradiction. Similarly, $u'\not\sim y$. But now, $\{u,v_3,u',x,y\}$ induces a $P_3\cup P_2$, again a contradiction. Thus we assume $y\in N_{\{3,4\}}$. Then $v_1$-$x$-$y$ is a $P_3$. It implies that $A$ is $P_2$-free, a contradiction to $\omega(A)\ge3$.

Next, we prove that $D$ is complete to $N_{\{1,2\}}\cup N_{\{3,4\}}$. Suppose not. By symmetry, there exist two nonadjacent vertices $x$ and $u$ such that $x\in N_{\{1,2\}}$ and $u\in D$. By Lemma \ref{nonadjcent}, we have that $x$ has a neighbor $y$ such that $y\not\sim u$. It is clear that $y\in N_{\{1,2\}}\cup N_{\{3,4\}}$. If $y\in N_{\{1,2\}}$, let $y'$ be a vertex in $A$ such that $y'\sim u$, then $\{y',u,v_3,x,y\}$ induces a $P_3\cup P_2$, a contradiction. So, $y\in N_{\{3,4\}}$. Let $x',y'$ be an edge in $G[A]$. Then $\{x,y,v_3,x',y'\}$ induces a $P_3\cup P_2$, a contradiction. Therefore, $D$ is complete to $N_{\{1,2\}}\cup N_{\{3,4\}}$.

Note that $\{v_1,v_2\}\cup N_{\{1,2\}}$ and $\{v_3,v_4\}\cup N_{\{3,4\}}$ are both cliques, and $N_{\{1,2\}}\cup N_{\{3,4\}}\cup \{v_1,v_2,v_3,v_4\}$ is complete to $N_{\{1,2,3,4\}}$ as $N_{\{1,2,3,4\}}=C\cup D$. Let $\omega(G[N_{\{1,2,3,4\}}])=\omega_0$. Since the vertices of an odd hole or an odd antihole cannot be divided into two cliques, by the famous strong perfect theorem, we have that $\chi(G[N_{\{1,2\}}\cup N_{\{3,4\}}\cup \{v_1,v_2,v_3,v_4\}])\le \omega(G[N_{\{1,2\}}\cup N_{\{3,4\}}\cup \{v_1,v_2,v_3,v_4\}])\le \omega(G)-\omega_0$. Therefore, $\chi(G)\le \chi(G[N_{\{1,2,3,4\}}])+\chi(G[N_{\{1,2\}}\cup N_{\{3,4\}}\cup \{v_1,v_2,v_3,v_4\}]\cup G[A])+\chi(G[B])\le 2\omega(G[N_{\{1,2,3,4\}}])+\omega(G)-\omega_0+3 \le 2\omega_0+\omega(G)-\omega_0+3$ as $A$ is anticomplete to $N_{\{1,2\}}\cup N_{\{3,4\}}\cup \{v_1,v_2,v_3,v_4\}$  and $\chi(G[B])\le3$.  Since $\omega(A)\ge3$, we have that $\omega(G)-\omega_0\ge3$, which implies that $\chi(G)\le 2\omega(G)$, a contradiction.

\medskip

{\bf Case 2} $1\le\omega(G[A])\le 2$.

\medskip

We divide the proof process into two cases depending on $N_{\{1,2\}}\cup N_{\{3,4\}}\neq\emptyset$ or $N_{\{1,2\}}\cup N_{\{3,4\}}=\emptyset$.

\medskip

{\em Subcase 2.1}  $N_{\{1,2\}}\cup N_{\{3,4\}}\neq\emptyset$.

\medskip

Suppose $\omega(G[A])= 2$. From the arguments in Case 1, the difference is that $A$ may not be complete to $N_{\{1,2,3,4\}}$. But now, since $N_{\{1,2\}}\cup N_{\{3,4\}}\neq\emptyset$, we have $\omega(\{v_1,v_2\}\cup N_{\{1,2\}})\geq 3$ or $\omega(\{v_3,v_4\}\cup N_{\{3,4\}})\geq3$. Then $\chi(G)\le \chi(G[N_{\{1,2,3,4\}}])+\chi(G[N_{\{1,2\}}\cup N_{\{3,4\}}\cup \{v_1,v_2,v_3,v_4\}]\cup G[A])+\chi(G[B])\le 2\omega_0+\omega(G)-\omega_0+3\le 2\omega(G)$, a contradiction. So, we may assume that $\omega(G[A])= 1$.

Suppose $\omega(G[C])\ge 2$. Recall that $C$ is complete to $N_{\{1,2\}} \setminus \{x_0\}$ and $N_{\{3,4\}} \setminus \{y_0\}$ for some $x_0 \in N_{\{1,2\}}  $, $y_0 \in N_{\{3,4\}} $. And $x_0,y_0$ are both anticomplete to $C$. Thus, $\omega(G[C\cup\{x_0,y_0\}])=\omega(G[C])$ as $\omega(G[C])\ge 2$. Therefore, $\chi(G)\le \chi(G[C\cup\{x_0,y_0\}])+\chi(G[(N_{\{1,2\}}\setminus \{x_0\})\cup (N_{\{3,4\}}\setminus \{y_0\})\cup \{v_1,v_2,v_3,v_4\}]\cup G[A])+\chi(G[B\cup D])\le 2\omega_1+\omega(G)-\omega_1+3\le 2\omega(G)$, a contradiction.

Suppose $\omega(G[C])=1$. $\chi(G)\le \chi(G[C\cup\{x_0\}])+\chi(G[(N_{\{1,2\}}\setminus \{x_0\})\cup (N_{\{3,4\}}\setminus \{y_0\})\cup \{v_1,v_2,v_3,v_4\}]\cup G[A])+\chi(G[B\cup D])+\chi(\{y_0\})\le 1+ \omega(G)-1+\omega(G)-1+1=2\omega(G)$, a contradiction.

\medskip

{\em Subcase 2.2}  $N_{\{1,2\}}\cup N_{\{3,4\}}=\emptyset$.

\medskip

Now, $A$ is anticomplete to $\{v_1,v_2,v_3,v_4\}$, $\chi(G)\le \chi(G[N_{\{1,2,3,4\}}])+\chi(G[\{v_1,v_2,v_3,v_4\}]\cup G[A])+\chi(G[B])\le 2\omega_0+ 2+\omega(G)-\omega_0$. Since $\omega(G)-\omega_0\ge2$, we have that $\chi(G)\leq 2\omega(G)$, a contradiction.

\medskip

{\bf Case 3} $A=\emptyset$.

\medskip

If $C$ is complete to $N_{\{1,2\}}\cup N_{\{3,4\}}$, then $\chi(G)\le\chi(G[C])+\chi(G[N_{\{1,2\}}\cup N_{\{3,4\}}\cup \{v_1,v_2,v_3,v_4\}])+\chi(G[B\cup D])\le 2\omega_0+2(\omega(G)-\omega_0)=2\omega(G)$, a contradiction. So, we may assume that $C$ is not complete to $N_{\{1,2\}}$ ($C$ is not complete to $N_{\{3,4\}}$). Let $x_0\in N_{\{1,2\}}$  ($y_0\in N_{\{3,4\}}$) such that $x_0$ is not complete to $C$  ($y_0$ is not complete to $C$). Then $x_0$ is anticomplete to $C$ ($y_0$ is anticomplete to $C$).

If $\omega_0\ge2$, then $\chi(G)\le \chi(G[C\cup\{x_0,y_0\}])+\chi(G[(N_{\{1,2\}}\setminus \{x_0\})\cup (N_{\{3,4\}}\setminus \{y_0\})\cup \{v_1,v_2,v_3,v_4\}])+\chi(G[B\cup D])\le 2\omega_0+2(\omega(G)-\omega_0)=2\omega(G)$, a contradiction.

If $\omega_0=1$, then $\chi(G)\le \chi(G[C\cup\{x_0\}])+\chi(G[(N_{\{1,2\}}\setminus \{x_0\})\cup (N_{\{3,4\}}\setminus \{y_0\})\cup \{v_1,v_2,v_3,v_4\}])+\chi(G[B\cup D])+\chi(\{y_0\})\le 1+ \omega(G)-1+\omega(G)-1+1= 2\omega(G)$, a contradiction.

This completes the proof of Lemma \ref{c7}.\qed

\medskip

\noindent{\em Proof of Theorem \ref{P3}}: By Lemma \ref{c7}, we may assume that $N_{\{1,2,3,4\}}$ is a clique with $|N_{\{1,2,3,4\}}|=\omega_0$. Now, we divide the proof process into two cases depending on $A\ne\emptyset$ or $A=\emptyset$.

\medskip

{\bf Case 1} $A\ne\emptyset$.

\medskip

First, we prove that $A$ is complete to $N_{\{1,2,3,4\}}$. Suppose not. Let $u\in N_{\{1,2,3,4\}}$ and $v\in A$ such that $u\not\sim v$. By Lemma \ref{nonadjcent}, it follows that there exists a vertex $v'\in N(v)$ and $v'\not\sim u$. By Lemma \ref{c0} and $A$ is anticomplete to $N_{\{1,2\}}\cup N_{\{3,4\}}$, we have that $v'\in N_{\{1,2,3,4\}}$, this contradicts that $N_{\{1,2,3,4\}}$ is a clique. So, $A$ is complete to $N_{\{1,2,3,4\}}$.

Second, we prove that $N_{\{1,2\}}$ and $N_{\{3,4\}}$ are both almost complete to $N_{\{1,2,3,4\}}$. By symmetry, we only need to prove $N_{\{1,2\}}$ is almost complete to $N_{\{1,2,3,4\}}$. Suppose that there exist two vertices in $N_{\{1,2\}}$ are not complete to $N_{\{1,2,3,4\}}$, say $x,y$. Since $N_{\{1,2\}}$ is a clique, $x\sim y$. Let $a\in A$. If $x$ and $y$ are both nonadjacent some vertex in $N_{\{1,2,3,4\}}$, say $u$, then $\{a,u,v_3,x,y\}$ induces a $P_3\cup P_2$ as $a$ is complete to $N_{\{1,2,3,4\}}$ and is anticomplete to $N_{\{1,2,\}}\cup N_{\{1,2,\}}$, a contradiction. So, $x$ and $y$ cannot be both nonadjacent some vertex in $N_{\{1,2,3,4\}}$. Then there exist two vertices $u_1$ and $u_2$ such that $x\not\sim u_1,x\sim u_2,y\not\sim u_2,y\sim u_1$, which implies that $\{x,y,u_1,u_2,v_4\}$ induces a house, a contradiction. So, $N_{\{1,2\}}$ is almost complete to $N_{\{1,2,3,4\}}$ and $N_{\{3,4\}}$ is almost complete to $N_{\{1,2,3,4\}}$.

Next, we prove that if $x\in N_{\{1,2\}}$ is not complete to $N_{\{1,2,3,4\}}$, then there exists a vertex $y\in N_{\{3,4\}}$ such that $y$ is not complete to $N_{\{1,2,3,4\}}$ and $x\sim y$, and $|A|=1$. Suppose $u\in N_{\{1,2,3,4\}}$ such that $x\not\sim u$. By Lemma \ref{nonadjcent}, there exists $x'\in N(x)$ such that $x'\not\sim u$. Since $A\cup B\cup (N_{\{1,2\}}\setminus\{x\})$ is complete to $N_{\{1,2,3,4\}}$, we have that $x'\in N_{\{3,4\}}$.  If $x'\sim u$, then $\{x',x,v_2,u,v_3\}$ induces a house, a contradiction. So, $x'$ is not anticomplte to $N_{\{1,2,3,4\}}$, that is to say, $(N_{\{3,4\}}\setminus\{y\})$ is complete to $N_{\{1,2,3,4\}}$. Consequently, suppose $A$ has two vertices, say $a,b$. Then $\{a,b,u,x,x'\}$ induces a $P_3\cup P_2$ if $a\not\sim b$ and $\{v_1,x,x',a,b\}$ induces a $P_3\cup P_2$ if $a\sim b$, both are contradictions. So, $|A|=1$.

Furthermore, it follows that $N_{\{1,2\}}$ is complete to $N_{\{1,2,3,4\}}$ if and only if $N_{\{3,4\}}$ is complete to $N_{\{1,2,3,4\}}$.

Suppose $N_{\{1,2\}}$ is complete to $N_{\{1,2,3,4\}}$. Then $\chi(G)\le \chi(G[N_{\{1,2,3,4\}}])+\chi(G[N_{\{1,2\}}\cup N_{\{3,4\}}\cup \{v_1,v_2,v_3,v_4\}\cup A])+\chi(G[B])\le \omega_0+\omega(G)-\omega_0+\omega(G)-\omega_0\leq 2\omega(G)$ as $G[N_{\{1,2\}}\cup N_{\{3,4\}}\cup \{v_1,v_2,v_3,v_4\}\cup A]$ is perfect.

Therefore, we may assume that $x\in N_{\{1,2\}}$ is not complete to $N_{\{1,2,3,4\}}$ and $y\in N_{\{3,4\}}$ is not complete to $N_{\{1,2,3,4\}}$, and $A=\{a\}$. Note that $N_{\{1,2,3,4\}}\ne\emptyset$, let $u\in N_{\{1,2,3,4\}}$. We can assign the color of $u$ to $x$ and a new color to $y$. Then $\chi(G)\le \chi(G[N_{\{1,2,3,4\}}])+\chi(G[(N_{\{1,2\}}\setminus\{x\})\cup (N_{\{3,4\}}\setminus\{y\})\cup \{v_1,v_2,v_3,v_4\}\cup A])+\chi(G[B])+1\le \omega_0+\omega(G)-\omega_0+\omega(G)-\omega_0+1\leq 2\omega(G)$ as $\omega_0\ge1$, a contradiction.

\medskip

{\bf Case 2} $A=\emptyset$.

\medskip

\pf By Lemma \ref{c0}, we have that $B$ is complete to $N_{\{1,2,3,4\}}$. Note that $\chi(G[B])=\omega(G[B])$. Thus $\chi(G[B\cup N_{\{1,2,3,4\}}])\leq \omega(G)$. And
$G[(N_{\{1,2\}} \cup (N_{\{3,4\}} \cup \{v_1,v_2,v_3,v_4\}]$ is perfect. Hence
$\chi(G)\le \chi(G[N_{\{1,2,3,4\}}\cup B])+\chi(G[N_{\{1,2\}}\cup N_{\{3,4\}}\cup \{v_1,v_2,v_3,v_4\}]) \le  \omega(G)+\omega(G)\leq 2\omega(G)$, a contradiction.

This completes the proof of Theorem \ref{P3}.\qed

\end{document}